\documentclass[11pt]{article}

\usepackage{fullpage}
\usepackage{amsfonts} 

\usepackage{amsthm}

\usepackage{enumitem}

\usepackage{graphics}

\usepackage{here} 

\usepackage{epsfig}
\usepackage{epstopdf}

\usepackage{tikz}
\usetikzlibrary{arrows}
\usetikzlibrary{shapes}
\usetikzlibrary{positioning}

\usepackage{xcolor}

\definecolor{darkbrown}{rgb}{.5,.1,.1} 

\usepackage{refcheck} 
\norefnames	
\nocitenames 	
\ignoreunlbld 
\setoffmsgs 

\def\ni{\noindent}
\def\bs{\bigskip}
\def\ms{\medskip}

\newtheorem{theorem}{Theorem} 
\newtheorem{proposition}[theorem]{Proposition}
\newtheorem{lemma}[theorem]{Lemma}

\newtheorem{remark}[theorem]{Remark}

\def\M{\overline M}
\def\bk{\backslash}
\def\s{\setminus}

\begin{document}

\title{\Large 
\large Short rewriting, and geometric explanations related to the active bijection, for:\\
 Extension-lifting bijections for
oriented matroids, \\
\large
by S. Backman, F. Santos, C.H. Yuen, arXiv:1904.03562v2 {\small (October 29, 2023)}}
\author{Emeric Gioan}

\date{\vspace{-1mm}\small \today}

\maketitle




\vspace{-5mm}
\noindent {\bf Abstract.} 
For an oriented matroid $M$, and given a generic single element extension and a generic single element lifting of $M$, the main result of \cite{extlift} provides a bijection between bases of $M$ and certain reorientations of $M$ induced by the extension-lifting.
%
This note is intended 
to somehow clarify and precise the geometric setting for this paper in terms of oriented matroid arrangements and oriented matroid programming,
to describe and prove the main bijective result in a short simple way, 
and to show  how it consists of combining two direct bijections and a central bijection, 
which is  
the same as 
a special case - practically uniform - of the bounded case of the active~bijection~\cite{GiLV04, GiLV09}.\break
(The relation with the active bijection is addressed in \cite{extlift} in an indirect and more complicated way.)
\bs

This note contains a short rewriting, in a few simple results, for the main construction of~\cite{extlift},
completed with illustrations, geometric explanations, and remarks on how it fits the setting of \cite{GiLV04, GiLV09}.
%
As an introduction, the reader can already have a look at Figure \ref{fig:bijection} and its caption, which sum up the whole constructions and bijections. 
For precise statements and explanations, we first recall below usual notions and geometric interpretations from oriented matroid theory. See \cite{OM} for details on this theory.
Concerning the relation of \cite{extlift}
with
the active bijection, the main reference is \cite{GiLV09}, which addresses the bounded case in general oriented matroids. The uniform case, which is very close to the setting used here, is addressed in \cite{GiLV04}. (Let us also mention that those constructions were first introduced in \cite{Gi02}, and that, beyond the very particular scope of this paper, the active bijection is the subject of a series of papers and constructions, a short sum up of which is given~in~\cite{Gi22}.) 
Theorem~\ref{thm:main} in the end recalls the main theorem of \cite{extlift} and also rephrases it in terms of the active~bijection.

For brevity, we can write for instance $B\cup f$ or $B\s f$ instead of $B\cup \{f\}$ or $B\s \{f\}$,~respectively.

Let $M$ be an oriented matroid on $E$. For the sake of geometric interpretations, an \emph{arrangement of $M$} is a signed pseudosphere arrangement representing $M$. As usual, covectors correspond to faces in the arrangement and cocircuits to vertices ($0-$dimensional faces).
For $A\subseteq E$, we denote by $-_AM$ the oriented matroid obtained from $M$ by reorienting $A$, and we consider that $M$ has $2^E$ such \emph{reorientations} $-_AM$, $A\subseteq E$, even if the resulting oriented matroids can be equal. In this way, acyclic reorientations of $M$ bijectively correspond to regions in an arrangement of $M$. 
In \cite{extlift}, a reorientation $-_AM$ is described as an ``orientation'' in $\{+,-\}^E$ where the $-$ signs form the set~$A$. 
In this note, we rather stick to standard terminology and formalism.
(Another useful translation is the following:
in \cite{extlift}, a signed subset is conformal with an orientation, when, to us, it yields a positive subset in the corresponding reorientation.)

For a basis $B$ of $M$ and $b\in B$, the \emph{fundamental cocircuit} of $B$ with respect to $b$, denoted by $C^*(B;b)$, is the cocircuit of $M$ whose support is the complementary of the closure of $B\setminus b$, and with a positive sign for $b$.
Geometrically, it corresponds to the vertex intersection of elements in $B\setminus b$, on the positive side of $b$.
For a basis $B$ and $e\not\in B$, the \emph{fundamental circuit} of $B$ with respect to $e$, denoted by $C(B;e)$, is the circuit of $M$ whose support is contained in $B\cup e$, and with a positive sign for $e$. Recall that $b\in C(B;e)$ if and only if $e\in C^*(B;b)$ (in terms of supports), and that, in this case, $C(B;e)\cap C^*(B;b)=\{e,b\}$, which implies by orthogonality that the sign of $b$ in $C(B;e)$ and the sign of $e$ in $C^*(B;b)$ are opposite.

An \emph{extension} of $M$ with respect to an element $f\not\in E$ is an oriented matroid $\M$ on $E\cup f$ such that $\M\bk f=M$.
It is called \emph{generic} if every circuit of $\M$ containing $f$ is spanning, that is, the circuits of $\M$ containing $f$ are the sets $B\cup f$ where $B$ is a basis of $M$. 
Geometrically, given an arrangement for $M$, a generic extension consists in forming an arrangement of the same dimenion with a new element $f$ in general position, that is, that does not meet any vertex ($0$-dimensional face).
For a generic extension, cocircuits $D$ of $M$ bijectively correspond to cocircuits $D\cup f$ of $\M$ containing $f$.
Then, for an unsigned cocircuit $D$ of $\underbar M$, we may denote by $\sigma^*(D)$ the signed cocircuit of $M$ whose corresponding cocircuit in $\underbar M$ is $D$ and whose corresponding cocircuit in $\M$ has a positive sign for~$f$.
The mapping $\sigma^*$, called \emph{signature} of the extension, determines the extension with respect to~$f$; geometrically, it determines the position of $f$ in the arrangement of $M$, by listing the vertices contained in its positive side.
(In \cite{extlift}, $\sigma$ serves as a notation in two different ways, we use the one of the main theorem statement.)
See \cite[Chapter 7]{OM} for more details on this extension~construction.


A \emph{lifting} of $M$ with respect to an element $g\not\in E$ is a dual notion (the dual of a lifting of $M$ is an extension of the dual of $M$). It is an oriented matroid $\M$ on $E\cup g$ such that $\M/ g=M$.
It is generic if cocircuits of $\M$ containing $g$ are the sets $(E\s B)\cup g$ where $B$ is a basis of $M$.
Geometrically, a generic lifting consists in forming an arrangement for $\M$ of one more dimension than $M$, such that $M$ is the arrangement induced on $g$, and such that each vertex of $\M$ not in $g$ is the intersection of exactly $rank(M)$ elements, forming a basis of $M$. In this case,  circuits $C$ of $M$ bijectively correspond to circuits $C\cup g$ of $\M$ containing $g$.
Then, for an unsigned circuit $C$ of $\underbar M$, we may denote by $\sigma(C)$ the signed circuit of $M$ whose corresponding circuit in $\underbar M$ is $C$ and whose corresponding circuit in $\M$ has a positive sign for $g$. It yields the \emph{signature} $\sigma$ of the lifting.
\ms


\vspace{-5mm}
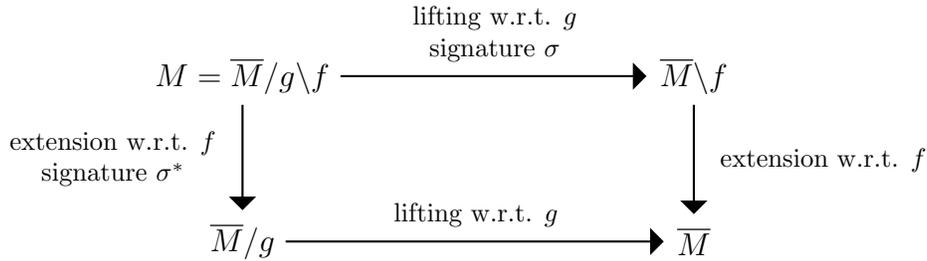
\begin{figure}[H]
\centering
\def\hdistance{6cm}
\def\vdistance{2.2cm}
\scalebox{1}{
\begin{tikzpicture}%
[
->,>=triangle 90, thick,shorten >=1pt,auto, node distance=\hdistance,  thick,  
main node/.style={rectangle,font=\sffamily\large}
   ]
   
  \node[main node] (1a) {$M=\M/g\bk f$};
  \node[main node] (1b) [right of=1a] {$\M\bk f$};
   \node[main node] (2a) [below of=1a, node distance=\vdistance] {$\M/g$};
  \node[main node] (2b) [right of=2a] {$\M$};

\path[every node/.style={}]
	
	(1a) edge [->] node [bend right] 
		{\small 
		\begin{tabular}{c}
			lifting w.r.t. $g$\\
			signature $\sigma$
		\end{tabular}} 
	(1b)

	(2a) edge [->] node [bend right] 
		{\small 
		\begin{tabular}{c}
			lifting w.r.t. $g$
		\end{tabular}}  
	(2b)

	(1a) edge [->] node [left] 
		{\small 
		\begin{tabular}{c}
			extension w.r.t. $f$\\
			signature $\sigma^*$
		\end{tabular}} 
	(2a)
	
	(1b) edge [->] node [right] 
		{\small 
		\begin{tabular}{c} 
				extension w.r.t. $f$
		\end{tabular}} 
	(2b)

	;		
\end{tikzpicture}
}
\vspace{-4mm}
\caption{Sum up of the extension-lifting construction of $\M$ from $M$.}
\label{fig:diagram}
\end{figure}
\vspace{-3mm}

From the usual notions above, the useful and practical setting for the construction addressed in this note is the following.
Given a generic extension  of $M$ with respect to $f$ (with signature $\sigma^*$), and a generic lifting  of $M$ with respect to $g$ (with signature $\sigma$), one can define an oriented matroid $\M$ on $E\cup\{f,g\}$ such that $\M/g\bk f=M$, which is both a generic extension of $\M\bk f$ with respect to $f$ and a generic lifting of $\M/g$ with respect to $g$. In this situation, we say that $\M$ is \emph{a generic extension-lifting of $M$} with respect to $(f,g)$ (and/or, equally, with respect to $(\sigma^*,\sigma)$). (Note that, by construction, $f$ and $g$ are not loops nor coloops of $\M$.)
See Figure \ref{fig:diagram} for a sum up of the involved oriented matroids. See the left part of Figure \ref{fig:bijection} for an illustration.

This setting allows us to handle the two initial extension and lifting of $M$ in a consistent way in a single $\M$, but let us observe that the choice of such an $\M$ is not unique (that is, $\M$ is not uniquely determined by $\sigma$ and $\sigma^*$). Indeed, starting with $M$ of rank $r$, and with the lifting of $M$ with respect to $g$ (given by $\sigma)$, yielding the uniquely determined $\M\bk f$ (of rank $r+1$), then any generic extension with respect to $f$ consistent with the extension $\M/g$ of $M$ with respect to $f$ (that is, consistent with $\sigma^*$) can be used. Geometrically: one can use any addition of an element $f$ in general position (in rank $r+1$) whose intersection with $g$ (in rank $r$) is fixed by the uniquely determined  $\M/g$ (given by $\sigma^*$). 

Amongst the possible choices for $\M$, let us say that the extension-lifting is \emph{compliant} if, geometricallly, there is no vertex in the positive side of $g$ and negative side of $f$ (this is always possible: in other words, $f$ is chosen to be as close to $g$ as possible, with prescribed intersection with $g$).

\def\hfig{7.9cm}
\begin{figure}[H]
\begin{center}
\hfil
\includegraphics[height=\hfig]{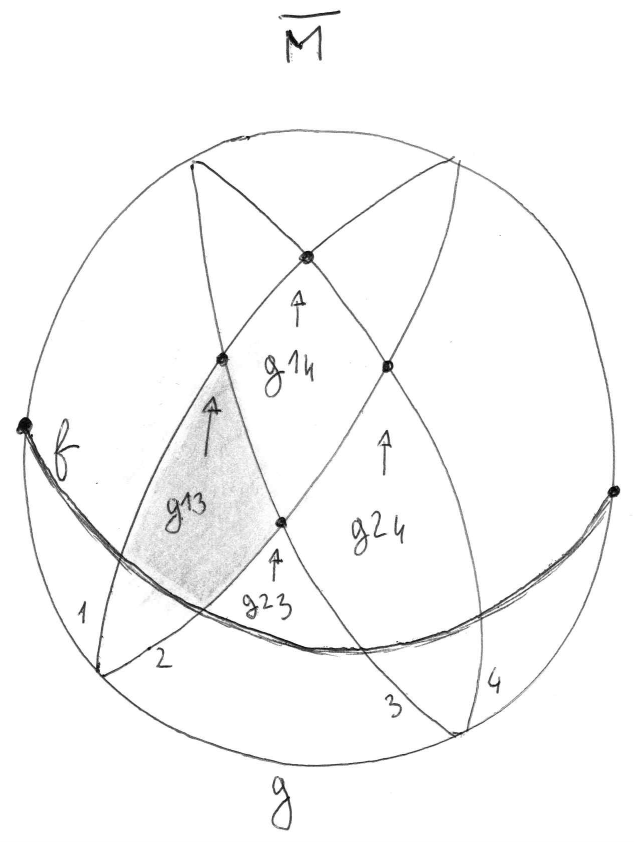}
\hfil
\includegraphics[height=\hfig]{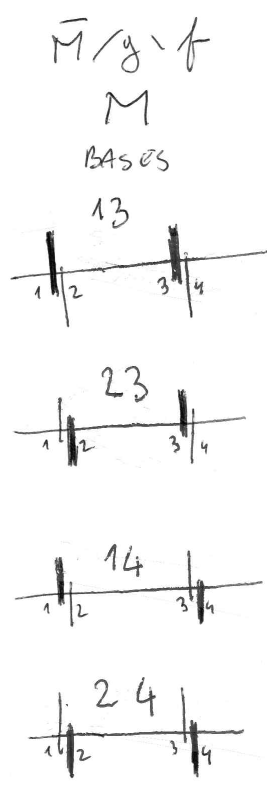}
\hfil
\includegraphics[height=\hfig]{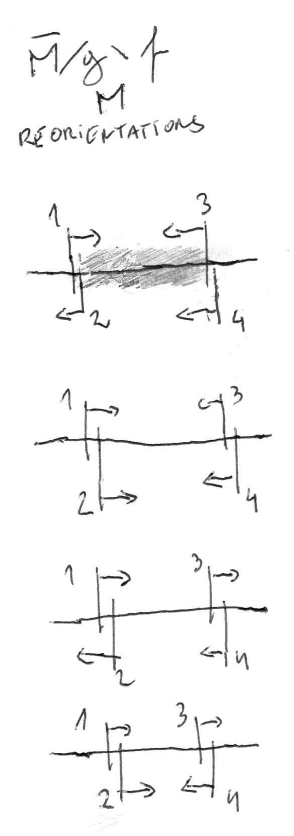}
\hfil
\end{center}
\vspace{-8mm}
\caption{
This figure sums up the whole constuctions and bijections. (1) On the left: an arrangement of $\M$ (the half-sphere given by the positive side of $g$, where the element $g$ is depicted as a circle), where $\M$ is a generic extension-lifting of $M$ with respect to $(f,g)$. One can see an arrangement of $M=\M/g\bk f$ as the arrangement induced on the element $g$ and removing $f$: $M$ consists of four elements on a line, with two parallel elements $1$ and $2$, and two parallel elements $3$ and $4$. (2) In the middle: bases of $M$. One can see that each basis corresponds to exactly one vertex in $\M$ given by the intersection point of its elements, and to the basis of $\M$ obtained by adding $g$. This is the first direct bijection, Proposition \ref{prop:bijbas}. (3) On the right: compatible reorientations of $M$, indicated with arrows pointing towards the positive side of the element. The grey one correwponds to the grey region in the left picture, given by the same positive sides in $\M$ plus a positive side for $f$ and $g$. One obtains in this way the bounded regions of $\M$ with respect to $(f,g)$. This is the second direct bijection, Proposition \ref{prop:compat}. Observe that in this arrangement, the position of $f$ is chosen so that $\M$ is compliant, meaning that bounded regions are all completely on the positive side of $f$ and correspond to regions not touching $g$. (4) Back on the left again: in each bounded region is written the basis associated by the central bijection of Proposition \ref{prop:bijcentral},
that is, its fully optimal basis in terms of the active bijection with respect to $E=g<f<E\s\{f,g\}$,
that is also, its optimal basis in terms of the oriented matroid program $(M,g,f)$, that is also, in ``common'' geometric terms, the basis whose associated intersection point is the ``farthest'' from $f$ in the direction of the ``infinity'' $g$. This bijection combined with the two previous ones yields the bijection from \cite{extlift} between bases and compatible reorientations~of~$M$  (as shown by the four rows in the middle and right parts).
}
\label{fig:bijection}
\end{figure}

%
%

\begin{proposition}[direct bijection between bases of $M$ and bases/vertices of $\M$]
\label{prop:bijbas}
Let $M$ be an oriented matroid on $E$ of rank $r$ and let $\M$ be a generic extension-lifting of $M$ with respect to~$(f,g)$.
%
\begin{enumerate}[label=(\roman*)]
\itemsep=-0.5mm
\partopsep=-0.5mm 
\topsep=0mm 
\parsep=0mm
\item \label{it1}
In an arrangement of $\M$, every vertex not in $g$ and not in $f$ is the intersection of exactly $r$ elements, forming a basis of $M$. 

\item \label{it2}
The mapping $B\mapsto (E\s B)\cup \{f,g\}$ is a bijection between bases of $M$ and cocircuits of  $\M$ containing $f$ and $g$ (they correspond to the vertices addressed in item \ref{it1}).
\item \label{it3} 
The mapping $B\mapsto B\cup g$ is a bijection between bases of $M$ and bases of $\M$ containing $g$ and not containing $f$ (their fundamental cocircuits for $g$ are the cocircuits addressed in item \ref{it2}).
\item \label{it4}
For a basis $B\cup g$ of $\M$ containing $g$ and not $f$,  $f\in C^*(B;b)$ for every $b\in B$, and $g\in C(B;e)$ for every $e\in E\s B$.
In particular, $C(B;f)$
has $B\cup \{f,g\}$ as support. (Geometrically, $f$ is in general position with respect to the vertices/faces intersections of elements in $B\cup g$.) 
\end{enumerate}
\end{proposition}

\begin{proof}
These assertions (closely related ot each other) all directly come from the fact that the lifting with respect to $g$ is generic and the extension with respect to $f$ is generic (as described above).
\end{proof}

\begin{remark}\rm
\label{rk:uniformoid}
Thanks to the properties above, we shall speak of a ``uniformoid'' situation, allowed by the generic extension and lifting. 
Indeed, these properties show that the vertices not in $g$ nor in $f$ are organized just like a uniform oriented matroid, in the sense that each vertex not in $g$ determines exactly one basis when $g$ is added, even if $\M/g$ is not uniform (and, furthermore, $f$ is in general position with respect to them).  Actually, if $\M$ is a uniform oriented matroid with $\M/g\bk f=M$ then the assertions hold in the exact same way. 
This is what will allow us to have a simple bijection in what follows, actually a very special case of the bounded case of the active~bijection.
\end{remark}

\begin{lemma}%
[same properties as in \cite{extlift} to define compatibility]
\label{lem:compat}
For an extension $\M$ of an oriented matroid $M$ with respect to $f$, and a reorientation $-_AM$, the two following properties are~equivalent.%
\vspace{-1mm}
\begin{itemize}
\itemsep=-1mm
\partopsep=-1mm 
\topsep=0mm 
\parsep=0mm
\item For every positive cocircuit $D$ of $-_AM$,
the corresponding cocircuit  $D\cup f$ in $-_A\M$ has a positive sign for $f$; 
\item  $-_f-_A\M$ is totally cyclic.
\end{itemize}
For a lifting $\M$ of an oriented matroid $M$  with respect to $g$, and a reorientation $-_AM$, the two following properties are equivalent.%
\vspace{-1mm}
\begin{itemize}
\itemsep=-1mm
\partopsep=-1mm 
\topsep=0mm 
\parsep=0mm
\item For every positive circuit $C$ of $-_AM$,
the corresponding circuit  $C\cup g$ in $-_A\M$ has a positive sign for $g$; 
\item  $-_g-_A\M$ is acyclic.
\end{itemize}
\end{lemma}

\begin{proof}
We prove the property for extensions, the property for liftings is a dual property.
Up to replacing $M$ with $-_AM$, we can assume that $A=\emptyset$.
Assume $-_f\M$ is totally cyclic. Then the signed set on $E$ with only $+$ signs and a $-$ sign for $f$ is a maximal vector of $\M$. Let $D$ be a positive cocircuit of $M$. Then, by orthogonality of vectors and cocircuits, the corresponding cocircuit  $D\cup f$ in $\M$ must have a positive sign for $f$.
Conversely, assume  $D\cup f$ is positive in $\M$ for every positive cocircuit $D$ of $M$. If $-_f\M$ is not totally cyclic then there exists a positive cocircuit $D$ in $-_f\M$. If $f\in D$ then $D\s f$ is a positive cocircuit of $\M\bk f$ and $f$ should be signed $-$ in $D$ which is a contradiction. If $f\not\in D$, then $D$ is not a cocircuit of $M$ since the extension is generic (geometrically, $D$ corresponds to the intersection of $f$ with a pseudoline in $\M$). Then there exists a cocircuit $D'$ in $\M$ with $f\in D'$, negative on $f$, and conformal with $D$ hence positive on $D'\s f$ ($D'$ is on this pseudoline, see \cite[Prop 7.1.4]{OM} for details) yielding the same contradiction.
\end{proof}


In the first situation addressed in Lemma \ref{lem:compat}, we say that $-_AM$ is \emph{compatible with the extension $\M$ with respect to $f$}. In the second situation, we say that  $-_AM$ is \emph{compatible with the lifting $\M$ with respect to $g$}.
Finally, when both an extension and a lifting are considered, if $-_AM$ is both compatible with the extension of $M$ with respect to $f$ with signature $\sigma^*$, and compatible with the lifting of $M$ with respect to $g$ with signature $\sigma$, then we say that 
 $-_AM$ is \emph{$(\sigma^*,\sigma)$-compatible},
 or also that
 $-_AM$ is \emph{compatible with the extension-lifting $\M$}, 
 for such an  extension-lifting $\M$ with respect to $(f,g)$, or equally with respect to $(\sigma^*,\sigma)$.

\ms

The reorientations of $\M$ addressed above can easily be geometrically described.
Let $\M$ be an oriented matroid on a set $E\cup \{f,g\}$, where $f,g\not\in E$, $f\not= g$, and $f$ and $g$ are not loops nor coloops.
Let us consider a reorientation $-_A\M$ such that $f,g\not\in A$, $-_A\M\bk f$ is acyclic and $-_A\M/g$ is totally cyclic. The direct geometric translation of this situation, in an arrangement of $\M$, is that $-_A\M$ corresponds to a region of $\M\bk f$ (on the positive side of $g$), possibly cut by $f$, and such that the vertices of the region belonging to $g$ are on the negative side of $f$.
In this situation, let us say that $-_A\M$ is a \emph{bounded region of $\M$ with respect to $(f,g)$}. 

\begin{proposition}[direct bijection between compatible reorientations of $M$ and bounded regions of $\M$]
\label{prop:compat}
Let $M$ be an oriented matroid and $\M$ be a generic extension-lifting of $M$ with respect to $(f,g)$.
The mapping $-_AM\mapsto -_{fg}-_A\M$ yields a bijection between reorientations of $M$ compatible with the extension-lifting $\M$, and bounded regions of $\M$ with respect to $(f,g)$.
\end{proposition}

\begin{proof}
Direct by the definitions and translation above.
\end{proof}

Continuing from the above setting, let us consider a bounded region $-_A\M$ with respect to $(f,g)$.
The following discussion is intended to explain how one can handle this setting in interrelated~ways.

First, this situation can be addressed in terms of oriented matroid programming. See \cite[Chapter 10]{OM} for details. 
In these terms, for a bounded region $-_A\M$ with respect fo $(f,g)$, the program $(-_A\M,g,f)$ has the property of being ``feasible and bounded''. (In the real case, or intuitively, it means that $\M$ delimits a proper region of the space, where the optimal vertices with respect to the objective function $f$ are at finite distance with respect to the element at infinity $g$.) Then, the fundamental result of oriented matroid programming applies: such a program has an \emph{optimal basis}, that is a basis $B$ such that $B$ contains $g$ but not $f$, $C^*(B;g)\setminus f$ is positive in $-_A\M$, and $C(B;f)\setminus g$ is positive in $-_A\M$. Geometrically, the vertex corresponding to  $C^*(B;g)$ on the positive side of $g$ is an \emph{optimal vertex} (in the real case alluded to above, it is an optimal vertex with respect to $f$ in the usual linear programming sense; in the oriented matroid case, it is a vertex with no outgoing edge in the skeleton of the region directed from the negative to the the positive side of $f$).

Second, the situation can be addressed in terms of the bounded case of the active bijection, as defined and studied in \cite {GiLV04, GiLV09}. We assume, for consistency of the settings, that $E$ is linearly ordered, with any ordering satisfying $E=g<f<E\s\{f,g,\}$. Furthermore, we assume that the extension-lifting is compliant (which can always be assumed and changes nothing to the final construction). Then, a bounded region $-_A\M$ with respect fo $(f,g)$ is ``$g$-bounded'' in the sense of \cite{GiLV09}, meaning it is a region of the arrangement of $\M$ such that none of its vertices belong to $g$, which is seen as an ``element at infinity''. (Equivalently,  in terms of ``activities'', the reorientation has activity $0$ and dual-activity $1$ with respect to the ordering of $E$, but these considerations are out of the scope of the present paper: in general, the active bijection allows to relate all bases and all reorientations in all ordered oriented matroids with respect to ''activities'', hence its name.)

In this setting, the main result of \cite{GiLV09} states that bounded regions of $\M$ with respect to $(f,g)$ are in bijection with their \emph{fully optimal bases}. (Formal definitions are given below in Remark \ref{rk:compare}.)
%
As noted in \cite[Definition 3.1]{GiLV09}, a fully optimal basis is, in particular, an optimal basis in the above sense. 
The crucial point for the present paper is that we handle a ``uniformoid'' situation, as described in Proposition \ref{prop:bijbas} and  Remark \ref{rk:uniformoid}, meaning, geometrically, that
each bounded region region has exactly one optimal vertex, and each optimal vertex determines exactly one optimal basis.
In this case, it is easy to see that the notion of fully optimal basis and optimal basis coincide.
Indeed, for an optimal basis, and for a fully optimal basis as well, the fundamental cocircuit of $g$ is an optimal vertex. Then, because of the  ``uniformoid'' situation, this vertex determines the basis, hence the two bases are the same. (In general, the notion of fully optimal basis is more complicated than that of optimal basis, as it involves the whole ordering of $E$ and all fundamental circuits and~cocircuits.)

Finally, the point is that, in this ``uniformoid'' situation, a bijection between bases and bounded regions can be simply handled in the same way as in the uniform case addressed in \cite{GiLV04}. 
%
Concretely, we get 
Proposition \ref{prop:bijcentral} below.
At the same time, it is a special (practically uniform) case of the bounded case of the active bijection \cite{GiLV09}, and it can be obtained by applying oriented matroid programming to each bounded region, as noted in \cite{extlift} and \cite{GiLV04} as well.
All this is illustrated on~Figure~\ref{fig:bijection}.%


\begin{remark}\label{rk:compare}\rm
\vspace{-1mm}
Let us add some more precise results to compare \cite{extlift} and \cite{GiLV04, GiLV09} in the above setting.
\vspace{-2mm}
\begin{enumerate}[label=(\roman*)]
\itemsep=-1mm
\partopsep=-1mm 
\topsep=0mm 
\parsep=0mm
\item\label{rk2:it1} 
According to \cite[Definition 3.1]{GiLV09}, a basis $B$ of an ordered oriented matroid $\M$ is fully optimal if and only if:
(a) for every $e\in E\s B$, the signs in $C(B; e)$ of $e$ and $\hbox{min}(C(B; e))$ are
opposite, and
(b) for every $b\in B\s g$, the signs in $C^*(B; b)$ of $b$ and $\hbox{min} (C^*(B; b))$ are
opposite.
In the present setting for $\M$, by Proposition \ref{prop:bijbas} \ref{it4} and by $E=g<f<\dots$, 
 the above definition is equivalent to: (a) for every $e\in E\s B$, the signs in $C(B; e)$ of $e$ and $g$ are
opposite, and
(b) for every $b\in B\s g$, the signs in $C^*(B; b)$ of $b$ and $f$ are
opposite.
One can observe that this formulation is very similar to the one used in \cite[Definition 1.1]{extlift}, recalled in Theorem \ref{thm:main} below. 
(One needs to reorient $f$ and/or $g$ to fit with the definitions of compatibility used in \cite{extlift}.)

\item\label{rk2:it2} 
Furthermore, since the support of $C^*(B;g)$ is $(E\s B)\cup\{g\}$ 
and the support of  $C(B; f)$ is $B\cup\{f\}$, 
the above definition is also equivalent to: (a)  $C^*(B;g)$ is positive, and
(b) $C(B;f)$ is positive except on $g$.
This 
expresses the fact that $B$ is an optimal basis as recalled above. It is the formulation used in \cite[Theorem 3.1]{extlift}, in \cite[Definition 3.1]{GiLV04}, and in Proposition \ref{prop:bijcentral}~below.

\item As one can see in Proposition \ref{prop:bijcentral}, in the above setting, 
(fully) optimal bases of bounded regions 
are
exactly 
bases containing $g$ but not $f$.
In general, in \cite{GiLV09}, bounded regions are in bijection with 
bases with internal activity $1$ and external activity $0$, also called ``uniactive internal bases'' (which necessarily contain $g$ but not $f$). Again, we do not need considerations on activities here.
Still, it is easy to check that bases containing $g$ but not $f$ exactly correspond to uniactive internal bases  here, again because of the ``uniformoid'' property of Proposition~\ref{prop:bijbas}~\ref{it4}: since every fundamental cocircuit contains $f$, then the only internally active elements is $g$, and since every fundamental circuit contains $g$, then there is no externally active~element.%

\item Starting from a basis that contains $g$ but not $f$, it is easy to build a bounded region for which this basis if (fully) optimal. One just has to reorient elements so that the definition is satisfied.
This is done in the same way in Proposition \ref{prop:bijcentral} below, in \cite{extlift}, and in \cite[Definition 3.1]{GiLV04}.
(This is done more generally for uniactive internal bases in \cite[Proposition 4.2]{GiLV09}.
The fact that it yields a bounded region is easy to check, and is stated in general 
in \cite[Proposition~3.2]{GiLV09}.)



\item The compliant case for the extension-lifting $\M$ is the case in which constructions in \cite{GiLV04, GiLV09}, in \cite{extlift}, and in oriented matroid programming exactly coincide.  Changing the position of $f$ to a non-compliant position  yields no change for the constructions in \cite{extlift} and in oriented matroid programming,
while it yields changes for the construction in \cite{GiLV04, GiLV09}, in which $f$ is part of the structure in a more precise way. 
However, 
starting from a non-compliant situation, 
one can always transform it into a compliant situation by changing the position of $f$, and 
directly retrieve the constructions for \cite{extlift} or oriented matroid programming of the initial situation.
(Also, for the final bijection of Theorem \ref{thm:main}, one can choose $\M$ to be in a compliant situation.)%
\end{enumerate}
\end{remark}


\vspace{-2mm}

\begin{proposition}[central bijection between bases/vertices of $\M$ and bounded regions of $\M$, equivalent to a special - practically uniform - case of the bounded case of the active bijection]
\label{prop:bijcentral}
Let $\M$ be a generic extension-lifting of an oriented matroid $M$ with respect to $(f,g)$.
\smallskip

For a bounded region $-_A\M$ of $\M$ with respect to $(f,g)$, there exists a unique basis $B$ of $\M$ containing $g$ but not $f$ such that $C^*(B;g)\setminus f$ is positive in $-_A\M$ and $C(B;f)\setminus g$ is positive in~$-_A\M$.
For a basis $B$ of $\M$ containing $g$ but not $f$, the reorientation $-_A\M$ such that $C^*(B;g)\setminus f$ is positive in $-_A\M$ and $C(B;f)\setminus g$ is positive in $-_A\M$ is a bounded region of $\M$  with respect~to~$(f,g)$.

We thus have in this way a bijection between bounded regions of $\M$ with respect to $(f,g)$ and bases of $\M$ containing $g$ but not $f$.
\smallskip

Furthermore, the basis $B$ associated with a bounded region $-_A\M$ is the unique optimal basis of the oriented matroid program $(M,g,f)$. When the extension-lifting $\M$ is compliant, 
this basis $B$ is the fully optimal basis of the bounded region $-_A\M$, and the above bijection is the same as the active bijection for bounded regions of $\M$,
with respect to any ordering $E=g<f<E\setminus\{f,g\}$.
(When 
it is not compliant, then one can consider a compliant extension-lifting $\M'$ where $f$ is replaced with $f'$ such that, geometrically, $f\cap g=f'\cap g$, and obtain the bijection for $\M$ by replacing $f'$ with $f$.)%
\end{proposition}

\vspace{-2mm}

\begin{proof}[Several proofs]
This bijection can be seen as a special case of the bounded case of the active bijection, that is, this proposition can be seen as a corollary of the constructions in \cite{GiLV09} and
\cite[Theorem 4.5]{GiLV09}. 
%
Formally, the translation in terms of the active bijection follows from the translations stated in Remark \ref{rk:compare}. 
(In the non-compliant situation, $f$ is replaced with $f'$ in a compliant position, 
which changes its sign in the involved signed subsets, positive in $C^*(B;g)$ in the compliant case, but not the bijection between regions of $\M\bk f$ and their associated bases.)

However, one does not need to use this ``big result'' of \cite{GiLV09} to prove the present result, as it is an easy special case, very close to the uniform case, as discussed previously. 
One can actually give various direct proofs of the uniqueness and existence of the basis $B$ (yielding the bijection),
as done  in \cite{GiLV04} and \cite{extlift}. Note that one does not necessarily need to use oriented matroid programming results either, even though the bijection can be naturally interpreted in those terms, as discussed~above.


{\sl --- Uniqueness.} 
A direct proof for the uniqueness part is done in \cite[Theorem 3.2, and Lemma 3.2.3]{GiLV04} and in \cite[Theorem 3.1]{extlift} in quite similar ways. It uses nothing but the definition of $B$ and some oriented matroid technique.
Let us give here another proof for the unicity 
based on oriented matroid programming.
Assume two bounded regions share the same optimal basis. Then they share the same optimal vertex $v$. Since this vertex is the intersection of exactly $rank(M)$ elements (cf. Proposition \ref{prop:bijbas}), then the two regions share a pseudoline in their border. In each region, the edge of the skeleton of the region contained in this pseudoline must be directed towards the vertex $v$, by properties of optimal bases. But on the other hand all edges on the pseudoline are directed from the negative side of $f$ towards its positive side, and $f$ does not contain $v$, so this is a contradiction.

{\sl --- Existence.} The existence part also has various proofs, directly coming from known results. 
One can use the main theorem of oriented matroid programming \cite[Theorem 10.1.13]{OM} stating that there exists an optimal basis (cf. setting discussed above). This is the argument used in \cite[Theorem 3.1]{extlift}. One can also use 
 the known result that the sets have equal cardinalities thanks to a classical result by Zaslavsky \cite{Za75}  generalized by Las Vergnas \cite{LV78}, so that the existence is implied by the unicity. This 
is the argument used 
in \cite[Theorem 3.2]{GiLV04} and more generally 
in~\cite[Theorem~4.5]{GiLV09}.
\end{proof}

\begin{theorem}[Rewriting  {\cite[Theorem A]{extlift}}, combining the bijections of Propositions \ref{prop:bijbas}-\ref{it3}, \ref{prop:compat}, and~\ref{prop:bijcentral}]
\label{thm:main}

 Let $M$ be an oriented matroid. Let $\sigma$ and $\sigma^*$ be the signature of a generic lifting and a generic extension of $M$, respectively.
 
For a basis $B$ of $M$, 
let $O(B)$ be  the reorientation 
of $M$
such that  each
$e\not\in B$ is reoriented~if~it is negative in $\sigma(C(B; e))$, and each $e\in B$ is reoriented if it is negative in $\sigma^*(C^*(B;e))$. 
Then~the~map%

\centerline{$B\mapsto O(B)$}

\ni is a bijection between the set of bases of $M$
and the set of $(\sigma^*,\sigma)$-compatible reorientations of $M$.

Equivalently, this bijection can be defined in the following terms (and written in the inverse~way).
Denoting by $\alpha$ the fully optimal basis of a bounded ordered oriented matroid as in \cite{GiLV09}, defining $\M$ as a compliant generic extension-lifting of $M$ with respect to $(f,g)$ given by the signatures $(\sigma^*,\sigma)$, and choosing any linear ordering for $E$ such that $E=g<f<E\s\{f,g\}$, 
the~map  

\centerline{
$-_AM\mapsto \alpha(-_{fg}-_A\M)\s g$
}

\ni is 
a bijection
between 
the set of $(\sigma^*,\sigma)$-compatible reorientations of $M$ and the set of bases of $M$.%
%
%
\end{theorem}

%
%

\begin{proof}
This bijection is the direct combination of the three bijections of Propositions \ref{prop:bijbas}-\ref{it3},~\ref{prop:compat},~and~\ref{prop:bijcentral}. Note that the definition of $O(B)$ here is not formulated as in Proposition \ref{prop:bijcentral}, but let us check that it is equivalent by Proposition \ref{prop:bijbas} \ref{it4} and by orthogonality of circuits and cocircuits in $\M$. (This is the same reformulation  as noted in Remark \ref{rk:compare} \ref{rk2:it1} and \ref{rk2:it2}.) Denote $-_AM$ for $O(B)$ and $\bar B$ for the basis $B\cup g$ in $\M$. 
For $e\not\in B$, $e$ and $g$ have opposite signs in $C^*(\bar B;g)$ if and only if they have the same sign in $C(B;e)$.
Hence, $C^*(\bar B;g)\bk f$ is positive in $-_g-_A\M$ and in $-_{fg}-_A\M$.
For $b\in B$, $e$ and $f$ have opposite signs in $C(B;f)$ if and only if they have the same sign in $C^*(\bar B;b)$.
Hence, $C(\bar B;f)\bk g$ is positive in $-_f-_A\M$  and in $-_{fg}-_A\M$.
So the definition of $\bar B$ fits with Proposition~\ref{prop:bijcentral}.
\end{proof}

\vspace{-1mm}
The previous results and comments contain and precise the main result 
and most intermediate 
results and comments from~\cite{extlift}, including those on the active bijection. 
%
In \cite{extlift}, 
complementary results are given, not addressed in the present note:
a reciprocal result on how bijections as in Theorem \ref{thm:main} come from extension-liftings; a result stating that 
if the extension-lifting is 
lexicographic with respect to a linear ordering then $(\sigma^*,\sigma)$-compatible reorientations yield representatives of ``activity classes'' addressed in the general active bijection (see \cite{Gi22});
relations between the bijection of Theorem \ref{thm:main} and works by Backman, Baker and Yuen, works by Ding, and  oriented matroid~triangulations.


\vspace{-1mm}
\bibliographystyle{plain}
\vspace{-4.5mm}

\small


\end{document}